\documentclass[12pt,leqno]{article}
\usepackage{amsmath}

\newcounter{conjecture}
\newcounter{remark}
\newcounter{corollary}

\newenvironment{corollary}{\medskip{\bf Corollary \thecorollary.}
\addtocounter{corollary}{1}\em}{\rm}
\newtheorem{theorem}{Theorem}
\newtheorem{lemma}{Lemma}

\newcommand {\rrr}[1]{(\ref{#1})}

\def \be{\begin{equation}}
\def \ee{\end{equation}}
\def \bt{\begin{theorem}}
\def \et{\end{theorem}}
\def \bc{\begin{corollary}}
\def \ec{\end{corollary}}
\def \bea{\begin{eqnarray}}
\def \eea{\end{eqnarray}}
\def \bas{\begin{eqnarray*}}
\def \eas{\end{eqnarray*}}
\def \bl{\begin{lemma}}
\def \el{\end{lemma}}

\def \vski{\vspace{12pt}}

\def \({\left(}
\def \){\right)}

\def \bc{\begin{center} }
\def \ec{\end{center} }
\def \bs{\begin{slide} }
\def \es{\end{slide} }

\def\square{{\vcenter{\vbox{\hrule height.3pt
         \hbox{\vrule width.3pt height5pt \kern5pt
            \vrule width.3pt}
         \hrule height.3pt}}}}

\newcounter{cccases}
\begin{document}

\title{Several consequences of adding or removing an edge from an electric network}

\author{{\bf Greg Markowsky}\\
{\it Department of Mathematics,}\\
{\it Monash University, Melbourne, Australia}\\greg.markowsky@monash.edu.au \\
and \\
{\bf Jos\'e Luis Palacios}\\
{\it Electrical and Computer Engineering Department,}\\{\it The University of New Mexico, Albuquerque, NM 87131, USA}\\
jpalacios@unm.edu\\
}



\maketitle

\begin{abstract}
In certain instances an electric network transforms in natural ways by the addition or removal of an edge. This can have interesting consequences for random walks, in light of the known relationships between electric resistance and random walks. We exhibit several instances in which this can be used to prove facts or simplify calculations. In particular, a new proof is given for the formula for the expected return time of a random walk on a graph. We also show how hitting times can be calculated in certain instances when a network differs from a highly symmetric one by one edge.
\end{abstract}

\vski


\vski


\vski






\section{Introduction}

Let $G$ be a finite, connected graph with $n$ vertices and $m$ edges. We will use the standard notation $\sim$ to denote adjacency in the graph, and for $z \in G$ let $deg(z)$ denote the degree of $z$. Let $X_j$ denote simple random walk on $G$; that is, $X_j$ is the Markov chain taking values in the vertex set of $G$ with transition probabilities given by

\begin{equation*}
P(X_{j+1} = z | X_j = y) = \left \{ \begin{array}{ll}
\frac{1}{deg(y)} & \qquad  \mbox{if } y \sim z  \\
0 & \qquad \mbox{otherwise } \;.
\end{array} \right.
\end{equation*}

As is shown in the classic reference \cite{doysne}, and expanded upon in countless other works, many probabilistic quantities associated with $X_j$ bear interpretations in the theory of electric networks; these include hitting times, cover times, mixing rates, and numbers of visits to vertices. The purpose of this note is to give several examples of situations where the addition or removal of a judiciously chosen edge can greatly simplify calculations or prove new results.

\vski 

In the next section we give a new proof of a known result, the calculation of the expected value of the time it takes to return to the same vertex one starts from. In the subsequent section we show how hitting times may be calculated when a network differs by one edge from a highly symmetric one, in particular a walk-regular graph.

\section{Return times}

Let $T_z = \inf\{j \geq 0 : X_j =z \}$, and let $T_z^+ = \inf\{j \geq 1 : X_j =z \}$. The object of interest for us is the expected return time, $E_zT^+_z$. The following elegant theorem is well known.

\begin{theorem} \label{t1}
\begin{equation*}
E_zT^+_z = \frac{2m}{deg(z)}.
\end{equation*}
\end{theorem}

This result admits a considerable generalization. Consider each edge $(y,z)$ as a wire in a circuit with a given conductance $C_{yz}$, which is a nonnegative number which measures how easily electricity (and the random walk) passes along the edge. For each vertex $z$ let $C_z = \sum_{y \sim z} C_{yz}$. Let $X_j$ now be the Markov chain taking values in the vertex set of $G$ with transition probabilities given by

\begin{equation*}
P(X_{j+1} = z | X_j = y) = \left \{ \begin{array}{ll}
\frac{C_{yz}}{C_y} & \qquad  \mbox{if } y \sim z  \\
0 & \qquad \mbox{otherwise } \;.
\end{array} \right.
\end{equation*}

This is the random walk induced by the electric network, and simple random walk corresponds to taking conductances of 1 (or any positive constant) across each edge. It should be noted that this construction is in fact quite general, since any reversible Markov chain can be realized as such an induced random walk (see \cite[Ch. 9]{levin}). Let $C = \sum_{y \in G} C_y$. We then have the following extension of Theorem \ref{t1} (which in fact applies to infinite graphs as well under the assumption that $C$ is finite).

\begin{theorem} \label{t2}

\begin{equation*}
E_zT^+_z = \frac{C}{C_z}.
\end{equation*}
\end{theorem}

The standard method of proving Theorem \ref{t2} is to appeal to a result from Markov chain theory, namely that an irreducible Markov chain with a stationary distribution $\pi$ satisfies $E_zT^+_z = 1/\pi_z$; and then simply verifying that $\pi_z = \frac{C_z}{C}$ is the stationary distribution for the chain $X_j$ (see \cite[Sec. 1.7]{norris}). On the other hand, researchers studying electric resistance have uncovered many identities and bounds on such quantities as hitting times, commute times, and cover times (\cite{comcov}); mixing times (\cite[Ch. 4]{aldfill}); and edge-cover times (\cite{george}). It is therefore natural to search for a derivation of Theorems \ref{t1} and \ref{t2} which makes more use of the principles which relate electric resistance to random walks, especially in light of the statement of the second theorem. We now present such a proof, naturally of the more general result, Theorem \ref{t2}.

\vski

{\bf Proof of Theorem \ref{t2}}. Fix $z \in G$, and construct a new graph $\tilde G$ which contains $G$ as a subgraph by adding a vertex $\tilde z$ to $G$ with a single edge connecting $\tilde z$ to $z$. Assign a conductance of 1 to the new edge. Let $\tilde X_m$ denote a random walk on $\tilde G$ induced by the conductances present (the original ones in $G$, together with the edge with unit conductance connecting $z$ and $\tilde z$). For $x,y \in \tilde G$, let $\tilde R_{x,y}$ be the effective resistance within $\tilde G$ between $x$ and $y$, and let $\tilde T_y = \inf\{j \geq 0 : \tilde X_j =y \}$. Let $\tilde C$ be twice the sum of the conductances across all edges in $\tilde G$; note that $\tilde C=C+2$. It is known (\cite[Cor. 11, Ch. 3]{aldfill}) that

\begin{equation*}
E_{\tilde z}\tilde T_z + E_z\tilde T_{\tilde z} = \tilde C \tilde R_{z,\tilde z}.
\end{equation*}

However, it is trivial that $E_{\tilde z}\tilde T_z = 1$ (the walk beginning at $\tilde z$ has no choice but to pass to $z$ at time 1), and it is equally trivial that $R_{z,\tilde z} = 1$ (there are no paths from $\tilde z$ to $z$ except for along the edge connecting them). Making the necessary substitutions yields

\begin{equation} \label{hai}
E_z\tilde T_{\tilde z} = C+1.
\end{equation}

Now, at time $\tilde T_{\tilde z} - 1$ the walk $\tilde X_m$ must necessarily reside at $z$. Furthermore, at each visit to $z$ the walk $\tilde X_m$ has probability $\frac{1}{C_z+1}$ of passing to $\tilde z$ at the next step, and between visits to $z$ the walk performs excursions within $G$, which will each take an average of $E_zT^+_z$ steps. The number of excursions within $G$ before $\tilde T_{\tilde z}$ is a Bernoulli trial with probability $\frac{1}{C_z+1}$, and as is well known the expected number of such trials until first success is $C_z+1$; however the number of excursions will in fact be one less than the number of visits to $z$, since the walk begins at $z$. It follows then that the expected number of excursions within $G$ before $\tilde T_{\tilde z}$ will be $C_z$. Adding $1$ to record the final step from $z$ to $\tilde z$ we obtain

\begin{equation} \label{shang}
E_z\tilde T_{\tilde z} = C_z E_zT^+_z+1.
\end{equation}

Equating the right-hand sides of \rrr{shang} and \rrr{hai} completes the proof of Theorem \ref{t2}. $\bullet$

\vski

{\bf Remark:} The key idea of attaching a new vertex to a vertex $z$ in a graph and starting a random walk there, armed with the knowledge that the first step must be to $z$, appears also in different contexts in \cite[Lemma 3.1]{northpal} and \cite[Ex. 10.4]{levin}.

\section{Networks which are nearly symmetric}

   
As was alluded to in the introduction, calculations of effective resistances (and therefore of such quantities as hitting and commute times) is in general a difficult problem, even for networks of modest size. In this context, symmetry properties of the graph in question are particularly prized, as problems often become tractable. Strong results exist for hitting times on vertex-transitive and edge-transitive graphs (see \cite{palrenom}), but what happens when the graph in question is merely close to being vertex-transitive or edge transitive? In this section we discuss how expected hitting times may still be deduced in some cases, and  we will use walk-regular graphs to demonstrate the method; particular instances in which the graph is vertex-transitive or distance-regular follow as corollaries.
 
\vspace{12pt}
 
 A graph is {\it walk-regular} if, for any $k\ge 2$ the number of closed paths of length $k$ starting at a given vertex $x$ is the same for any choice of $x$.  Walk-regular graphs are clearly regular and the class or walk-regular graphs contains several well-known families of symmetric graphs: vertex-transitive, distance-regular and regular edge-transitive graphs (see \cite{godsilcomment} for the latter). For our purposes, the most salient feature  of a walk-regular graph $G$ (see \cite{georgeunpub}) is that given $a, b\in G$
\begin{equation*}
E_aT_b=E_bT_a.
\end{equation*}
This relation, together with Chandra et al.'s result (see \cite{comcov}) stating that
\begin{equation}
\label{chandra}
E_aT_b+E_bT_a=2|E|R_{ab},
\end{equation}
allows us to conclude that for a walk regular graph we have
\begin{equation}
\label{walk}
E_aT_b=|E|R_{ab}.
\end{equation}
 The Kirchhoff index of the simple connected undirected graph, defined as
\begin{equation*}
R(G)=\sum_{a<b}R_{ab}
\end{equation*}
has been proposed in several contexts as a measure of the robustness of a complex network.  Specifically, it has been proposed to study the value of $R(G)$ compared to that of $R(G^{'})$, where $G^{'}$ is the original graph minus one edge (see \cite{clemente2019novel}, \cite{wang2014improving}), and in all these contexts it is relevant to know the value of the increase from the original effective resistance $R_{ab}$ to $R^{'}_{ab}$, the effective resistance between $a$ and $b$ in $G^{'}$.  Due to Rayleigh's monotonicity law (see \cite{doysne}), removing an edge from the graph $G$ does not decrease the effective resistance between any two vertices $a, b$, i.e.,  $R^{'}_{ab} \ge R_{ab}$, but in general it is not known exactly what that increment is. However, when $a$ and $b$ are neighbors, if $ab$ is not a cut-edge, that is, if the removal of $ab$ still leaves the rest of the graph connected, we can prove the following
\begin{theorem} \label{resis}
For any graph $G$ such that $ab$ is not a cut-edge we have
\begin{equation*}
R^{'}_{ab}=\frac{R_{ab}}{1-R_{ab}}.
\end{equation*} 
\end{theorem}
{\bf Proof.} Consider the edge $ab$ and the rest of the graph as two resistors connected in parallel at $a$ and $b$.  By the rules of effective resistance of resistors in parallel we have
$$R_{ab}=\frac{1}{1+\frac{1}{R^{'}_{ab}}}.$$
 Solving for $R^{'}_{ab}$ ends the proof~~$\bullet$
 \vskip .1 in
 
Theorem \ref{resis} gives us the following corollaries.

\begin{corollary}
For any graph $G$ such that $ab$ is not a cut-edge we have
$$R^{'}_{ab}-R_{ab}=\frac{R_{ab}^2}{1-R_{ab}}.$$
\end{corollary}

\begin{corollary} \label{ssb}
Over all connected graphs on $n$ vertices with unit resistance across each edge, the difference $R^{'}_{ab}-R_{ab}$ attains a maximal value $\frac{(n-1)^2}{n}$ for the $n$-cycle $C_n$ and minimal value $\frac{4}{n(n-2)}$ for the complete graph $K_n$.
\end{corollary}

\vspace{12pt}

{\bf Proof}.  The function $f(x)=\frac{x^2}{1-x}$ is increasing in the interval $(0, 1)$, so we look for the largest and smallest effective resistance in an edge of a connected graph on $n$ vertices which is not rendered disconnected by the removal of an edge. It is evident from Rayleigh's monotonicity law that the smallest possible resistance is $\frac{2}{n}$ for the complete graph $K_n$.  It is less obvious that the largest possible value of $\frac{n-1}{n}$ occurs for the $n$-cycle $C_n$, but we may argue as follows. If $G$ has $n$ vertices and $ab$ is not a cut-edge then $ab$ must be contained in a shortest cycle of length $n' \leq n$, and if $G$ is not $C_n$ then $n'<n$. But then $R_{ab}$ is bounded below by $\frac{n'-1}{n'}$, which is the resistance across an edge in the cycle graph $C_{n'}$, and this is smaller than $\frac{n-1}{n}$, the corresponding value for $C_n$. The result follows. ~~$\bullet$
 \vskip .2 in
 
We note that, in the random walk setting, Corollary 2 applies only to simple random walk. We will henceforth work only with simple random walk; that is, we will assume that any graph considered has unit resistance across each edge. We now state the main result of this section.

\begin{theorem} \label{remove}
Suppose that $G$ is a walk-regular graph with $|V|\ge 3$. Form a new graph $G^{'}$ by deleting an edge $(a,b)$ which is not a cut-edge. Then the expected hitting time increases from
$$E_aT_b=|E|R_{ab}$$
to
$$E_aT^{G^{'}}_b=\frac{(|E|-1)R_{ab}}{1-R_{ab}}=\frac{E_aT_b-R_{ab}}{1-R_{ab}}.$$
\end{theorem}   
{\bf Proof}. The fact that $E_aT_b=|E|R_{ab}$ is (\ref{walk}) above.  Now, when computing  $E_aT_b^{G^{'}}$ as the sum of all possible lengths of paths from $a$ to $b$ times the probabilities of their occurrences, one notices that deleting the edge $ab$ implies deleting all paths that arrive at $b$ in one jump from $a$.  But this is mirrored in the computation of $E_bT_a^{G^{'}}$,  as all paths arriving at $a$ in one jump from $b$ are deleted, and the probability with which that last jump occurs is the same, $\frac{1}{d}$ in each direction, because the graph is $d$-regular, so both computations of $E_aT_b^{G^{'}}$ and $E_bT_a^{G^{'}}$ are affected in the same way, and all other paths in the computation of the original hitting times are kept for the hitting times with the edge removed. Therefore 
$E_aT_b^{G^{'}}=E_bT_a^{G^{'}}$ and we can apply (\ref{chandra}) with the number of edges reduced to $|E|-1$ and $R_{ab}$ replaced by $R^{'}_{ab}$~~$\bullet$
\vskip .2 in
Let us discuss some examples. If $G$ is regular and edge-transitive, it was shown in \cite{palrenom} that for every edge $ab$ we have $E_aT_b=|V|-1$. (This implies, by (\ref{walk}), that for an edge $ab$ we have $R_{ab}=\frac{|V|-1}{|E|}$, a fact shown in \cite{foster1949average} for the smaller class of arc-transitive graphs). On the other hand, if we remove edge $ab$ then the new hitting time is
$$E_aT_b^{G^{'}}=(|E|-1)R^{'}_{ab}=\frac{(|E|-1)(|V|-1)}{|E|-|V|+1}.$$

In particular, in the $n$-cycle, for neighboring $a$ and $b$, $E_aT^G_b=n-1$ and $E_aT^{G^{'}}_b$ increases to $(n-1)^2$. In the complete graph $K_n$, $E_aT^G_b=n-1$ increases to $\displaystyle E_aT^{G^{'}}_b=\frac{n(n-1)-2}{n-2}$. In the 3-dimensional cube, $E_aT_b^G=7$ and $E_aT_b^{G^{'}}=15.4$.  In the {\it unitary Cayley graph} $Cay(Z_n, U_n)$ of the integers modulo $n$, $Z_n$ (and $U_n$ its set of units) where $(x,y)\in E$ if and only if $x-y$ is relatively prime with $n$,  the graph is $\phi(n)$-regular, where $\phi(n)$ is the Euler totient function and if $(a,b)\in E$ we have $E_aT_b=n-1$ whereas if we delete the edge $(a,b)$ we get
$$E_aT_b^{G^{'}}=\frac{(n-1)\left(\frac{n\phi(n)}{2}-1\right)}{\left(\frac{n\phi(n)}{2}-n+1\right)}.$$

\section{Acknowledgements}

We'd like to thank Tim Garoni for helpful conversations, as well as an anonymous referee for useful suggestions. The first author is grateful for support from Australian Research Council Grants DP0988483 and DE140101201.

\bibliographystyle{alpha}
\bibliography{prob}

\end{document}